\documentclass{article}

\begin{document}

\title{K\"ahler-Einstein metrics and algebraic structures on limit spaces}
\author{Simon Donaldson}
\date{\today}
\maketitle



\section{Introduction}
\newcommand{\bC}{{\bf C}}
\newcommand{\bR}{{\bf R}}
\newcommand{\bP}{{\bf P}}
\newcommand{\dbd}{\partial\overline{\partial}}
\newtheorem{thm}{Theorem}
\newcommand{\tomega}{\tilde{\omega}}

This is an expository article following the lines of the author's lecture at the Journal of Differential Geometry conference in Harvard, September 2014, and largely based on joint work with S. Sun. In particular, we will give an overview of the recent paper \cite{kn:DS2}.

Our theme is the formation of singularities in complex algebraic geometry and complex differential geometry. One motivation for this comes from existence questions such as Yau's conjecture relating the existence of K\"ahler-Einstein metrics on Fano manifolds to algebro-geometric \lq\lq stability''. This conjecture was proved by Chen, Donaldson and Sun \cite{kn:CDS} using the method of deformation of metrics with cone singularities along divisors. Recently, another proof has been given by Datar and Sz\'ekelyhidi \cite{kn:DSz} using the  continuity method. In a third direction, Chen and Wang have obtained strong results about K\"ahler-Ricci flow \cite{kn:CW}, and it seems likely that these will  lead to another proof of Yau's conjecture. In all these approaches, a central issue is that one wants to take a limit of a family of K\"ahler metrics and to endow this limit with an algebro-geometric meaning. This is the essential bridge between the differential  geometry and PDE discussion and the algebraic geometry.

Another motivation comes from moduli questions. There can be fundamental difficulties in forming a moduli space of complex structures on a given manifold, because in its natural topology the space of complex structures modulo diffeomorphisms need not be Hausdorff, but these difficulties do not appear for moduli spaces of Riemannian metrics. Suppose we have some moduli space ${\cal N}$ of K\"ahler-Einstein structures on a manifold, modulo diffeomorphism. Then (as we will recall in more detail below) in suitable cases the metric geometry may give us a topological compactification $\overline{{\cal N}}$ of ${\cal N}$. It is then natural to ask for an algebro-geometric interpretation of $\overline{{\cal N}}$ as a moduli space of varieties, with singularities allowed. A more detailed question is to ask how the  algebro-geometric  singularities are reflected in the behaviour of the  metrics.

To fix ideas, suppose that we have a family $X_{t}$ of hypersurfaces of degree $d$ in $\bC\bP^{n+1}$ parametrised by a variable $t$ in the unit disc in $\bC$, developing an ordinary double point singularity at $t=0$. Thus, in affine co-ordinates centred at the singular point,  $X_{t}$ is defined by a polynomial equation 
$P_{t}(z_{0}, \dots, z_{n})= 0$ where
$$  P_{0}(z_{0}, \dots, z_{n}) = \sum z_{i}^{2} + {\rm higher\  order}, $$
and without much loss of generality we could suppose that $P_{t}= P_{0}-t$. Suppose we know that for $t\neq 0$ the manifold $X_{t}$ has a K\"ahler-Einstein metric. What can we say about the behaviour of these metrics as $t\rightarrow 0$? When the dimension $n$ is $1$ or $2$ this is well-understood.

\begin{itemize}
\item {\it Dimension 1}\  We suppose $d\geq 4$ so the Euler characteristic   of the plane curves $X_{t}$  is negative. For $t\neq 0$ each $X_{t}$ has a metric of constant Gauss curvature $-1$. The behaviour as $t\rightarrow 0$ is modelled on the family of metrics on the topological cylinder obtained by dividing the upper half plane $\{ {\rm Im} w>0\}$ by the ${\bf Z}$-action generated by $w\mapsto \lambda w$ where $\lambda = e^{\vert t \vert}$. The metrics develop a long thin \lq\lq neck'' and in the limit can be thought of as splitting into a pair of hyperbolic cusps, corresponding to the fact that $X_{0}$ is locally reducible, with two branches passing through the singularity.

More generally, for $g\geq 2$ we have the famous Deligne-Mumford compactification $\overline{{\cal M}_{g}}$ of the moduli space ${\cal M}_{g}$ of Riemann surfaces of genus $g$ which can be constructed either algebro-geometrically via stable curves, or using hyperbolic surfaces with cusps \cite{kn:W}.

\item {\it Dimension 2}\  Suppose again that the degree $d$ of the surfaces is  at least $4$. By the existence theorem of Yau from the 1970's, for $t\neq 0$ the $X_{t}$ have  K\"ahler-Einstein metrics. If $d=4$ the Ricci curvature is $0$ and if $d\geq 5$ it is $-1$, say. (In fact the discussion extends to the case $d=3$, with positive Ricci curvature.) The metric picture starts from the fact that the model singularity $\{ \sum z_{i}^{2}=0\}\subset \bC^{3}$ can be identified with the quotient $\bC^{2}/\pm 1$. There is a Ricci-flat {\it Eguchi-Hanson metric} on the affine quadric $\{\sum z_{i}^{2}=1\}\subset \bC^{3}$ which is asymptotic to $\bC^{2}/\pm 1$. Scaling this gives a family of metrics on the quadrics $\{ \sum z_{i}^{2}=t\}$ and this furnishes the local model for the K\"ahler-Einstein metrics on $X_{t}$. The limit of these metrics is a K\"ahler-Einstein  orbifold metric on $X_{0}$.  

There are generalisations of this picture, involving  quotient singularities $\bC^{2}/\Gamma$ for $\Gamma\subset U(2)$ and Ricci flat \lq\lq gravitational instantons'' asymptotic to $\bC^{2}/\Gamma$. By results from the 1980's of Anderson, Nakajima and Tian a sequence of K\"ahler-Einstein surfaces with fixed volume and scalar curvature and with bounded diameter has a subsequence converging to an orbifold with a finite number of singular points of this type. 

\end{itemize}

The problems start in dimension 3 or higher. Consider the case when $n=3$ and, for definiteness, that the degree $d$ is $5$ (but the discussion applies for other degrees, certainly $d\geq 5$).  Thus we are considering a family of quintic Calabi-Yau 3-folds developing a double point singularity. There is a candidate for a local model, much as in dimension $2$. The model singularity $\{ \sum z_{i}^{2}=0\}\subset \bC^{4}$ is the affine cone (in the algebro-geometric sense) over the corresponding projective quadric surface $Q\subset \bC\bP^{3}$ which has a standard K\"ahler-Einstein metric of positive Ricci curvature. If $S$  is any K\"ahler-Einstein manifold with positive Ricci curvature and $Y_{S}\rightarrow S$ is a circle bundle corresponding to some  fractional power $K_{S}^{1/p}$ of the canonical bundle there is an induced  Sasaki-Einstein metric on $Y_{S}$. This is equivalent to a Ricci-flat K\"ahler  metric on the Riemannian cone $C(Y_{S})$. In the case at hand (taking the bundle corresponding to $K_{Q}^{-1/2}$) this yields a description of the affine quadric as a metric cone $C(Y_{Q})$   with Ricci-flat metric. Further, there is an explicit Ricci-flat {\it Stenzel metric} on the smooth affine quadric $\{\sum z_{i}^{2}=1\}\subset \bC^{4}$ which is asymptotic to the cone $C(Y_{Q})$. Scaling this gives a family of metrics on the quadrics $\{\sum z_{i}^{2}=t\}\subset \bC^{4}$ converging to the cone metric on the singular quadric as $t\rightarrow 0$. It is natural to expect that this gives the model for the local behaviour of the metrics on $X_{t}$  but at present it seems that there is no known proof that is the case.  (It a pleasure here for the author to thank  Mark Haskins, who has emphasised the central nature of this problem for some years  past.) 

Approaching the problem from a slightly different direction, it is known by work of Eyssidieux, Guedj and Zeriahi \cite{kn:EGZ} that there is a Calabi-Yau metric $\omega_{0}$ on the singular variety $X_{0}$ (with an appropriate definition of a metric on a singular variety). Results of Ruan and Zhang \cite{kn:RuZ} and  Rong and Zhang \cite{kn:RZ} show that the metrics $\omega_{t}$ converge to $\omega_{0}$ away from the singular point. The essential gap, at present,  is  the question whether the metric $\omega_{0}$ is modelled on $C(Y_{Q})$ near the singular point. More precisely, in terms that we will recall below, the question is whether the {\it tangent cone} of $(X_{0},\omega_{0})$ at the singular point is $C(Y_{Q})$.

We will now explain the crucial difference, from the differential geometric and PDE point of view, between the cases $n=2$ and $n\geq 3$. When $n=2$ the singular space $X_{0}$ is an orbifold and we can choose a reference orbifold metric $\tomega$. Then Yau's existence proof for K\"ahler-Einstein metrics on  manifolds extends with no essential change to the orbifold case and gives an orbifold K\"ahler-Einstein metric $\omega_{0}$ \cite{kn:Ko}. With this in place, gluing techniques can be used to construct K\"ahler-Einstein metrics on $X_{t}$ for small $t$, as in the work of Spotti \cite{kn:Sp}. That is, one starts with an approximate solution, made by gluing a scaled version of the Eguchi-Hanson metric to $(X_{0},\omega_{0})$, and uses an implicit function theorem to produce a small deformation which is  K\"ahler-Einstein. This leads to the desired description of the family of metrics $\omega_{t}$.

To try the same approach in dimension $3$, we can certainly choose a reference metric $\tomega$ in $X_{0}$, modelled on the cone $C(Y_{Q})$ near the singular point, and set up a continuity method for a $1$-parameter family $\tomega + i \dbd \phi_{s}$. The  difficulty comes in extending the $C^{2}$ estimate for $\phi_{s}$. Yau's approach uses a maximum principle argument which requires  a lower bound on the holomorphic sectional curvature of the reference metric $\tomega$. There is a variant, applying the maximum principle to the  Chern-Lu inequality, which  requires an upper bound for the curvature of the reference metric; see \cite{kn:JMR}, Section 7,  for example.  But in our situation the curvature of the cone $C(Y_{Q})$ is not bounded above or below, so neither argument is available. This is the fundamental difficulty. The existence proof of \cite{kn:EGZ} for the K\"ahler-Einstein metric on $X_{0}$ follows different lines,  and does not give precise information about the behaviour of the solution near the singularity.

 \section{Review of convergence theory and results for compact K\"ahler-Einstein manifolds}

If $A$ and $B$ are compact metric spaces the Gromov-Hausdorff distance $d_{GH}(A,B)$ ca be defined by saying that $d_{GH}(A,B)\leq \delta$ is there is a metric on the disjoint union $A\sqcup B$ extending the given metrics on $A,B$ and such that each of $A,B$ are $\delta$-dense. This leads to a notion of a sequence of compact metric spaces $A_{i}$ with Gromov-Hausdorff limit a metric space $A_{\infty}$.
More generally, if $A_{i}$ are not necessarily compact metric spaces and  if there are base points $a_{i}\in A_{i}$   the based spaces $(A_{i}, a_{i})$ converge to a based limit $(A_{\infty}, a_{\infty})$ if for each $R>0$ the metric balls $B(a_{i}, R)$ converge to $B(a_{\infty}, R)$ in the sense above. If the diameters of the $A_{i}$ have a  fixed upper bound then this reduces to the previous notion and the base points are irrelevant.

 A version of the fundamental compactness theorem of Gromov  states that if $A_{i}$ are complete Riemannian manifolds of fixed dimension, with a fixed lower bound on the Ricci curvature and with the usual induced metric space structure, and if $a_{i}\in A_{i}$ are any base points, then some subsequence has a based Gromov-Hausdorff limit. In particular this applies to sequences of Einstein manifolds with scalar curvature bounded below. Within this general framework a variety of things can happen. For one example, if we consider a family of Riemann surfaces $X_{t}$ with constant curvature $-1$ developing a long neck,  as in the previous section, and if we take base points in the middle of the neck then the based limit is the $1$-dimensional manifold $(\bR,0)$. Another kind of example occurs for the Gross-Wilson metrics on elliptically fibred $K3$ surfaces \cite{kn:GW} with the volume of the fibres shrinking to zero. If the K3 surfaces are scaled to have diameter $1$ then the Gromov-Hausdorff limit is the 2-sphere, with some metric structure.    In this article we want to avoid considering such \lq\lq collapsing'' phenomena. This can be achieved by requiring a definite lower bound \begin{equation}{\rm Vol}(B(a_{i},
r)\geq \kappa\end{equation} for some fixed $r,\kappa>0$ and all $i$, which we assume from now on. With this assumption, if  $A_{i}$ are $m$-dimensional Einstein manifolds, with constant Ricci curvature $\lambda$ for $\vert\lambda \vert \leq 1$ say, then there is an open dense \lq\lq regular set'' $A_{\infty}^{\rm reg}\subset A_{\infty}$ which is a smooth $m$-dimensional Einstein manifold and the convergence over compact subsets of $A_{\infty}^{{\rm reg}}$ can be realised by the standard notion of $C^{\infty}$ convergence of metric tensors. The complement $A_{\infty}\setminus A_{\infty}^{{\rm reg}}$ is a closed subset of codimension $4$ or more, in the Hausdorff sense.

Suppose that $A_{\infty}$ is a non-collapsed limit as above and $q$ is a point in $A_{\infty}$. For any sequence of real numbers $\mu_{j}\rightarrow \infty$ we consider the based metric spaces $(A_{\infty}, \mu_{j} d_{A_{\infty}}, q)$, i.e. we are just rescaling the metric by factors $\mu_{j}$. Passing to a subsequence , we can suppose this has a Gromov-Hausdorff limit---it follows from the definitions that this is also a based limit of  a rescaled subsequence of the original manifolds  $(A_{i'}, \mu_{i'}d_{A_{i'}}, b_{i'})$ for suitable base points $b_{i'}\in A_{i'}$. A basic theorem of Cheeger and Colding \cite{kn:CC} states that this limit is a metric cone $C(L)$, a {\it tangent cone} of $A_{\infty}$ at $q$. Uniqueness of these tangent cones, assuming only a lower bound on the Ricci curvature,  can fail \cite{kn:CC},  and uniqueness is not known in general even when the $A_{i}$ are Einstein.

With this background in place, we return to the complex geometry setting and recall a form of the main result of \cite{kn:DS1}. For fixed $D,n$ we consider the class ${\cal K}(n,D)$ of triples $(X,L)$ where $L\rightarrow X$ is a Hermitian holomorphic line bundle over a compact $n$-dimensional complex manifold and $h$ is a metric on $L$ with curvature form $-i\omega$ where $\omega$ is a K\"ahler-Einstein metric on $X$, with Ricci form $\lambda \omega$ where $\lambda$ is $1,0$ or $-1$ and the diameter of $(X,\omega)$ is bounded by $D$. If we have a sequence $(X_{i}, L_{i})$ in ${\cal K}(n,D)$ then by Gromov's result there is a Gromov-Hausdorff convergent subsequence of $X_{i}$ with the metric structures induced from $\omega_{i}$. So let us suppose that the original sequence converges in this sense to a limiting metric space $Z$. Since the volume of the $X_{i}$ is at least $(2\pi)^{n}$ the Bishop comparison theorem gives the non-collapsing bound (1).

 The essential fact, roughly stated,  is that $Z$ has the structure of a complex  variety and we give a number of statements which make this more precise. To formulate these it is useful to introduce the disjoint union
$$  {\cal X}= \left(\bigsqcup X_{i} \right) \sqcup Z, $$
which carries a natural topology making it a family ${\cal X}\rightarrow \Lambda$ over the topological space $$ \Lambda = \{1,1/2,1/3,\dots\} \cup \{0\}.$$

{\bf Complex analytic space} For an open set $U\subset Z$ we let ${\cal S}(U)$ be the set of functions on $U$ which extend to continuous functions on a neighbourhood of $U$ in ${\cal X}$, holomorphic in the $X_{i}$ variables. This is a pre-sheaf and we let ${\cal O}_{Z}$ be the corresponding sheaf. The first statement is that, after possibly passing to a subsequence, the ringed space  $(Z,{\cal O}_{Z})$ is a normal complex analytic space. The singular sets of $Z$ as a complex analytic space and as a Riemannian limit space coincide. 
(To see why a subsequence might be necessary, notice that we could take the \lq\lq opposite'' complex structure on, say, the $X_{i}$ with $i$ odd.)

Given the general Riemannian convergence results, it is an easy fact that the regular set $Z^{{\rm reg}}\subset Z$ has a complex structure.   One can also define ${\cal O}_{Z}$ as the push-forward of the sheaf of holomorphic functions on $Z^{{\rm reg}}$ under the inclusion.

{\bf Projective embeddings} 
 The statement is that (perhaps after passing to a subsequence) there  are fixed $k, N$ and a continuous map $\Phi: {\cal X}\rightarrow \bC\bP^{N}$ which on each $X_{i}\subset {\cal X}$ is an embedding defined by sections of $L_{i}^{k}\rightarrow X_{i}$ and restricts to a homeomorphism from $Z\subset {\cal X}$ to a normal projective variety $Z_{{\rm alg}}\subset \bC\bP^{N}$.  This means that the images $\Phi(X_{i})$ converge to $Z_{{\rm alg}}$ in the usual sense of complex algebraic varieties. In fact we get maps from $\Lambda$ to an appropriate Chow variety and Hilbert scheme.

{\bf Co-ordinate ring}

As above it is an easy fact that (possibly taking a subsequence) there is a limiting holomorphic line bundle $L_{\infty}$ over the regular set $Z^{{\rm reg}}\subset Z$. Then we have a ring
    $$ R= \bigoplus_{j=0}^{\infty} H^{0}(Z^{{\rm reg}}, L_{\infty}^{j}). $$
    The key statement is that $R$ is finitely generated and given this we can define $Z_{{\rm alg}}$ as a projective scheme
$ Z_{{\rm alg}}= {\rm Proj}(R)$. 

(This formulation makes use of results  of  C. Li \cite{kn:CL}.)

\

Of course all these statements are compatible in natural ways, and we have not attempted to formulate here a single statement which captures everything.

\section{Based limits and tangent cones}

In this section we discuss more recent results from \cite{kn:DS2}. Suppose that we have a sequence $(X_{i}, L_{i})$ in   ${\cal C}(n,D)$ as above and that we have base points $p_{i}\in X_{i}$ and scale factors $\mu_{i}\rightarrow \infty$. We scale the K\"ahler metrics $\omega_{i}$ by factors $\mu_{i}^{2}$ and taking a subsequence we can suppose that there is a based limit $(Z,p)$. The definition of a structure sheaf ${\cal O}_{Z}$ goes just as before and we have:
\begin{thm}
 The ringed space $(Z,{\cal O}_{Z})$ is a normal complex analytic space.
 \end{thm}
 
 The remaining statements will be given in an order which makes the exposition simpler, but which is different from the order in which the proofs go. 
 We define a ring $R$ whose elements are holomorphic functions on $Z$ of polynomial growth with respect to the distance from the base point $p$.
\begin{thm}
The ring $R$ is finitely generated and ${\rm Spec}(R)$ is a normal affine algebraic variety with underlying complex analytic space  $(Z,{\cal O}_{Z})$.
\end{thm}

These are the analogs of the first and third statements in the previous situation. We expect that  an analogue of the second statement should hold, so that in particular $Z$ is locally smoothable, but there are some technical difficulties in proving this.

We now specialise to our central concern which is the case of tangent cones. Let $Z$ be a limit  of K\"ahler-Einstein manifolds, either of bounded diameter (as considered in Section 2) or after scaling by $\mu_{i}\rightarrow \infty$, as considered above, and let $q$ be a point of $Z$.
\begin{thm}
  There is a unique tangent cone $C(Y)$ of $Z$ at $q$.
  \end{thm}

This tangent cone $C(Y)$ is itself a based and rescaled limit of spaces in ${\cal C}(n,D)$ so Theorem  2 applies to show that it has the structure of   an affine algebraic variety which we denote by $C(Y)_{{\rm alg}}$. From the Riemannian geometry point of view there is a Ricci-flat K\"ahler metric on the open, dense, regular set in $C(Y)$. Extending the uniqueness theorem of Berndtsson \cite{kn:BB} to the case of cones, we show that this metric is uniquely determined by the algebraic structure $C(Y)_{{\rm alg}}$. The essential question is
{\em can we determine $C(Y)_{{\rm alg}}$ by an entirely algebro-geometric procedure from the germ of $Z$ at $q$?} For example in the case we discussed in Section 1, when $Z$ has an ordinary double point at $p$, such a procedure would surely tell us that the tangent cone is $C(Y_{Q})$ and solve the problem we raised there.

In \cite{kn:DS2} we do not answer this essential question but we show that
$C(Y)_{{\rm alg}}$ is obtained from $(Z,q)$ by an algebro-geometric procedure which perhaps depends on certain extra data. The most important part of this data is a filtration of the local ring ${\cal O}_{q}$ of $Z$ at $q$. 

\begin{enumerate} \item  Suppose that we have a map 
$$  D: {\cal O}_{q}\rightarrow [0,\infty], $$
with \begin{itemize}
\item $D(f)=\infty$ if and only if $f=0$,
\item $D(f)=0$ if and only if $f(q)\neq 0$,
\item $D(f+g)\geq \max(D(f), D(g))$,
\item $D(fg)\geq D(f)+ D(g)$.
\end{itemize}

\item  Suppose that the image of ${\cal O}_{q}\setminus \{0\}$ is a discrete subset
of $[0,\infty)$, which we  order as $0<\nu_{1}<\nu_{2}<\nu_{3}\dots$. Then we define a filtration by ideals
$$     {\cal O}_{p}= {\cal F}_{0}\supset {\cal F}_{1}\supset {\cal F}_{2}\supset \dots , $$
by saying that $f\in {\cal F}_{i}$ if $D(f)\geq \nu_{i}$. Let $R_{D}$ be the associated graded ring
$$     R_{D}= \bigoplus_{j\geq 0} \frac{{\cal F}_{j}}{{\cal F}_{j+1}}. $$
We have an action of $(\bC,+)$ on $R_{D}$ with $s\in \bC$ acting as $\exp(\nu_{j}s)$ on $\frac{{\cal F}_{j}}{{\cal F}_{j+1}}$. 
 \item  Suppose that $R_{D}$ is finitely generated and each quotient $\frac{{\cal F}_{j}}{{\cal F}_{j+1}}$ is
finite-dimensional. 
\end{enumerate}

\

Then  $W={\rm Spec}(R_{D})$ is an affine scheme. The action of $\bC$ on the ring defines an action on $W$. There is an embedding $W\subset \bC^{N}$, a linear action of a complex torus
$T= (\bC^{*})^{m}$ on $\bC^{N}$ preserving $W$ and an element $v\in {\rm Lie}(T)$ such that the $\bC$ action on $W$ is given by $\exp(s v)$. Further  there is a neighbourhood $\Omega$ of $p$ in $Z$ and an embedding $\Omega\subset \bC^{N}$ such that $\exp(sv) \Omega$ converges to $W$ as $s$ tends to $+\infty$ through real values.

\

Using the metric structure on $Z$ we define
$$   D_{\rm metric}(f) = \limsup_{r\rightarrow 0} (\log r)^{-1} \max_{z\in B(p,r)} \log \vert f(z)\vert ,$$
interpreted as $+\infty$ if $f$ is identically zero. 
(In fact the limsup here can be replaced by lim: there is a well-defined \lq\lq order of vanishing'' at $p$ of holomorphic functions, with respect to the metric structure.) 
Then we have
\begin{thm} The map $D_{\rm metric}$ satisfies (1), (2), (3) above and $W={\rm Spec}( R_{D_{\rm metric}})$  is a normal affine variety.  The tangent cone $C(Y)_{{\rm alg}}$ is the central fibre of a $\bC^{*}$-equivariant degeneration of $W$, through affine varieties in $\bC^{N}$ invariant under the  $T$ action,  and the restriction of the $\bR$ action $\exp(sv), s\in \bR$ to $C(Y)_{{\rm alg}}$ corresponds to radial dilations on the metric cone. 
\end{thm}

The question then is whether, for the same $(Z,q)$, there could be  other such maps $D$ and degenerations of ${\rm Spec}(R_{D})$ to normal $T$-invariant varieties which admit Ricci-flat cone metrics. 

\

{\bf Remarks} 

\begin{itemize}\item There is a standard notion of a \lq\lq tangent cone'' in algebraic geometry. This corresponds to the filtration of the local ring by powers of the maximal ideal or to the scalar action $\exp(s)$ on $\bC^{N}$. In some cases, such as in the expected picture for ordinary double points, the tangent cone we are describing will coincide with this, but  in general it could be different.
\item The statement of Theorem 4 allows the possibility that the degeneration from $W$ to $C(Y)_{{\rm alg}}$ is trivial, so $C(Y)_{{\rm alg}}=W$. For example this is what we expect in the case of ordinary double points. 
\item The equivariant degenerations from $W$ to $C(Y)_{{\rm alg}}$ that we consider are analogs for cones of the degenerations (or \lq\lq test configurations'') appearing in the definition of K-stability for projective varieties. One could hope for a more algebro-geometric formulation in which $W$ is a semi-stable cone, in a suitable sense, and $C(Y)_{{\rm alg}}$ is the associated stable cone. 
\item An example where we expect that there is a non-trivial degeneration from $W$ to $C(Y)_{{\rm alg}}$ is a 3-fold singularity  with local equation
 $P(z_{1}, z_{2}, z_{3}, z_{4})=0$ where
 $$  P(\underline{z})= z_{1}^{2}+ z_{2}^{2}+ z_{3}^{2} + z_{4}^{4} + {\rm higher\ order\ terms}. $$
Here \lq\lq higher order'' means with respect to the weight $(2,2,2,1)$ which defines the $\bR$ action on $\bC^{4}$ and the corresponding $W$ is the weighted homogeneous singularity $\{z_{1}^{2}+ z_{2}^{2}+ z_{3}^{2} + z_{4}^{4}=0\}\subset \bC^{4}$. This does not admit a Ricci-flat cone metric, by results of Gauntlett, Martinelli, Sparks and Yau \cite{kn:GMSY}. However there is a degeneration through a family of equations
$$    z_{1}^{2}+ z_{2}^{2}+ z_{3}^{2} + \epsilon^{4} z_{4}^{4}=0, $$
to the singularity with equation $z_{1}^{2}+ z_{2}^{2}+ z_{3}^{2}=0$. This is the product of $\bC$ with the 2-dimensional  singularity $\bC^{2}/\pm 1$, which does admit a Ricci-flat cone metric. (See the more extended discussion in \cite{kn:DS2}, Section 3.3.) 
\item The whole discussion applies to the case of a smooth point $q$ in $Z$. In this case we do know by PDE arguments that the metric is smooth and the tangent cone is $\bC^{n}$. It is interesting to try to prove this algebro-geometrically using Theorem 4.  
\end{itemize}



\end{document}